\documentclass{amsart}
\usepackage{ifxetex}
\ifxetex
\usepackage{fontspec}
\usepackage{hyperref}
\else
\usepackage[utf8]{inputenc}
\usepackage[unicode,pdftex]{hyperref}
\fi

\usepackage{amsmath}
\usepackage{amsaddr}
\usepackage{amsthm}

\usepackage[style=numeric-comp,minnames=8,maxnames=99,isbn=false,url=false]{biblatex}
\usepackage[margin=1in]{geometry}

\usepackage{tikz}
\usetikzlibrary{%
  matrix,%
  calc,%
  arrows,%
  positioning%
}

\usepackage{bbold}
\usepackage{color}
\usepackage{url}
\usepackage{wrapfig}

\usepackage{algorithmic}
\usepackage{algorithm}

\newtheorem{theorem}{Theorem}

\newtheorem{proposition}[theorem]{Proposition}
\newtheorem{definition}[theorem]{Definition}

\newcommand\into\hookrightarrow
\newcommand\onto\twoheadrightarrow
\newcommand\from\leftarrow
\newcommand\infrom\hookleftarrow
\newcommand\onfrom\twoheadleftarrow

\newcommand\kk{\mathbb{k}}
\newcommand\N{\mathbb{N}}
\newcommand\R{\mathbb{R}}
\newcommand\eR{\bar{\mathbb R}}
\newcommand\X{\mathbb{X}}
\newcommand\Y{\mathbb{Y}}
\newcommand\Z{\mathbb{Z}}
\newcommand\mcF{\mathcal{F}}
\newcommand\mcG{\mathcal{G}}

\newcommand\mcI{\mathcal{I}}

\DeclareMathOperator\Dgm{\mathrm{Dgm}}
\DeclareMathOperator\Bij{\mathrm{Bij}}
\DeclareMathOperator\Amp{\mathrm{Amp}}
\DeclareMathOperator\Lip{\mathrm{Lip}}
\DeclareMathOperator\Pers{\mathrm{Pers}}
\DeclareMathOperator\Vect{\mathrm{Vect}}
\DeclareMathOperator\vC{\text{\v C}}
\DeclareMathOperator\VR{V\!R}

\newcommand\dfn[1]{\textbf{#1}}

\begin{document}

\title[Sketches of a platypus: a survey of persistent homology and its algebraic foundations]
{Sketches of a platypus: a survey of\\ persistent homology and its algebraic foundations}
\author{Mikael Vejdemo-Johansson}
\email{mvj@kth.se}
\address{Computer Vision and Active Perception Laboratory\\ KTH Royal
Institute of Technology\\ Teknikringen 14\\ 10044 Stockholm\\ Sweden \\ ~ \\
Artificial Intelligence Laboratory \\ Jo{\v z}ef Stefan Institute \\ Ljubljana, Slovenia}

\begin{abstract}
  The subject of persistent homology has vitalized applications of algebraic topology to point cloud data and to application fields far outside the realm of pure mathematics. The area has seen several fundamentally important results that are rooted in choosing a particular algebraic foundational theory to describe persistent homology, and applying results from that theory to prove useful and important results.

  In this survey paper, we shall examine the various choices in use, and what they allow us to prove. We shall also discuss the inherent differences between the choices people use, and speculate on potential directions of research to resolve these differences.
\end{abstract}

\maketitle

\textcite{johnstone_sketches_2002} named his book on topos theory ``Sketches of an elephant'', referencing a joke: three blind wise men encounter an elephant. They each try to describe it to each other. The wise man who caught hold of the elephant's trunk says ``An elephant is like a snake.''; the wise man holding the ear says ``An elephant is like a palm leaf.''; and the wise man holding its leg says ``An elephant is like a tree.''.

The joke is highly relevant to topos theory; which has its roots in logic, in geometry, and in topology, with the three perspectives being fundamentally different and enriching each other in surprising ways.

The title of this paper is similar, but different. The platypus is well-known to be a hybrid of an animal: sharing traits both with the phylum of birds and with the phylum of mammals. The field of persistent homology is in a similar situation to the platypus: there are at least two different viewpoints of what persistent homology should be, and they interact in sometimes unexpected ways.

\newpage
\tableofcontents

\newpage
\section{Introduction}
\label{sec:introduction}

Persistent homology is a technique that has sparked the birth of a new field of research; while several introductory texts have been written \cite{ghrist_barcodes_2008,eh-ct-09,carlsson_topology_2009,zomorodian_topology_2005}, and several good survey articles have been published \cite{edelsbrunner2008persistent,chazal_structure_2012}, most if not all target computer scientists or data scientists with an interest in topology.

The development of persistent homology and topological data analysis has been driven by algorithm development. In this paper, we will try to describe the field and its development with a view towards the different foundational viewpoints that have been leveraged to prove increasingly valuable results in the field.

As an alternative, this article proposes to be an introductory survey targeting mathematicians with an interest in the applicability, and with a specific view towards the applications of algebra in persistent homology. To our knowledge, there is one other article with a similar focus; the AMS Notices article by \textcite{w-wph-11}.

We assume that the reader is comfortable with the homology functor, basic category theory and homological algebra including the idea of an abelian category, and basic analysis including the idea of a Lipschitz function.

For the remainder of the paper, we will go through the various viewpoints and their strengths in order chosen for. To help the reader keep the descriptions in context, the article starts, here, with a very short overview of the upcoming contents.

\subsection{Foundations in use}
\label{sec:foundations-use}

There are two main genres of foundations in use, two cultures of ``persistent homology''. 

\begin{description}
\item[Filtered spaces] Persistent homology is about the effect of applying the homology functor to a filtration of topological spaces. Invariants describing the resulting homology diagrams help us construct tools for visualization and data analysis eventually allowing for the inference of topological structure for point clouds using specific constructions of filtered complexes that encode properties of point clouds.
\item[Representations of the reals] Persistent homology is about studying sublevel sets of real-valued functions on topological spaces. Such sublevel sets have -- for nice enough functions and spaces -- discretizations that allow us to adapt descriptions of finite diagrams of vector spaces to efficient descriptors. In particular, by using the ``distance from a set'' family of functions we can support inference of topological structures from point clouds.
\end{description}

Both these choices come with built in benefits as well as drawbacks. They give rise to different generalizations of the fundamental inference problem for point clouds sampled from a topological space, and they support different further constructions and proofs.

In particular, among the results that emerge from the two viewpoints, we will be discussing a selection in this paper:
\begin{description}
\item[Stability] The \textbf{representations of the reals} viewpoint allows us to prove a Lipschitz-style property for the inference process underlying the theory: there is a metric, the \dfn{bottleneck metric}, on the invariants of the diagrams of homology groups such that the distance between the homologies of the sublevel sets of two different functions is bounded by the $L_\infty$-distance between the functions. Evolutions in the exact definitions used for persistence lead to increasingly generous assumptions in this bound.
\item[Sub- and super- and iso-level sets] By modifying the constructions
used, we can get new constructions that allow us to study sequences of
super-level sets, of iso-level sets (or level-sets), and of the result
from collapsing sub- or super-level sets to a single point. In
particular, this brings us \dfn{extended persistence}, where no infinite
length intervals occur, and a number of topological features comes into
play, including Poincar\'e duality. Current technologies for iso-level sets tend to rely on zig-zag persistent homology (see below).
\item[Graded modules] The kinds of diagrams emerging from the \textbf{filtered spaces} viewpoint have the structure of graded modules over the polynomial ring $\kk[t]$. This recognition sparked both new algorithms for computing persistent homology with far less assumptions on the chosen coefficient ring, and a number of extensions of the fundamental constructions that we will mention below.
\item[Relax the \emph{filtration} requirement] In a seminal paper, \textcite{gabriel_unzerlegbare_1972} proved that the \emph{tameness} of the representation theory of quiver algebras depends only on the corresponding Dynkin diagram, not on the particular orientation of arrows in the quiver. Re-interpreting the diagrams of vector spaces emerging from the \textbf{filtered spaces} viewpoint as modules over quiver algebras rather than modules over $\kk[t]$ allows for inclusion maps that go both forwards and backwards producing \dfn{zig-zag persistent homology}, which has allowed for both a topological approach to statistical bootstrapping and concrete approaches to iso-level set persistent homology.
\item[More directions] The work by \textcite{carlsson_local_2007} on the topology of configurations of pixels in natural images relied on being able to vary several independent variables in the construction of the intermediate simplicial complexes studied. This inspired \textcite{carlsson_theory_2009} to study how these \emph{multi-dimensional} approaches can be handled. A straight generalization from graded $\kk[t]$-modules directs us to study modules over $\kk[t_1,t_2,\dots,t_n]$, which brings a whole range of theoretical and computational problems. Nevertheless, recent research seems promising. \cite{patriarca_presentation_2012,chacholski_combinatorial_2012}
\end{description}

Several results make specific reference to geometric complex constructions that are in common use in persistent homology. Since the choices of algebraic foundations seldom influence these constructions specifically, our description shall be brief and summary. Each of them requires a point cloud $L$ -- a finite subset of some metric space.
\begin{description}
\item[\v Cech complex] The \v Cech complex is an $\varepsilon$-parametrized simplicial complex defined as the nerve complex of the family of open $\varepsilon$-balls around the points of $L$. We write $\vC_\varepsilon(L)$ for the resulting filtered and parametrized simplicial complex.
\item[Vietoris-Rips complex] The Vietoris-Rips complex is the most widely used construction -- it is less dependent on dimension constraints than the $\alpha$-complex, and less computationally intensive to work with than the \v Cech complex. The Vietoris-Rips complex $\VR_\varepsilon(L)$ at $\varepsilon$ contains a simplex $(\ell_0,\dots,\ell_k)$ if for all $0\leq i\leq j\leq k$, $d(\ell_i,\ell_j)\leq\varepsilon$.
\item[$\alpha$-complex] The $\alpha$-complex is a powerful and very nice tool -- the intersection of the \v Cech complex and the Delaunay complex on a point cloud, it comes with strong theoretical guarantees. However, the computational complexity of the Delaunay complex means that the $\alpha$-complex is mainly of use in 2 and 3 ambient dimensions. $\alpha$-complexes were introduced by \textcite{edelsbrunner_shape_1983} in 2 dimensions and by \textcite{edelsbrunner_three-dimensional_1992} in 3 dimensions. The study of their Betti numbers by \textcite{robins_computational_2002} is one of the immediate precursors to the definition of persistent homology.
\item[Witness complexes] Witness complexes were introduced by \textcite{de_silva_topological_2004} as one approach to deal with the computational complexity of persistent homology. The construction uses a relatively small vertex set $L$ and an often far bigger \emph{witness set $W$}. Given a $k$-simplex $\sigma$ with vertices from $L$ and a points $w\in W$ we say that $w$ is an \dfn{$\alpha$-witness} of $\sigma$ if the vertices of $\sigma$ are all within $d_k(w)+\alpha$ of $w$, where $d_k(w)$ is the distance from $w$ and its $(k+1)$th nearest neighbour in $L$. We write $W_\alpha(L,W)$ for this simplicial complex. Witness complexes have been further studied by \textcite{chazal_towards_2008} and by \textcite{chazal_persistence_2012}.
\end{description}


\section{Persistence barcodes and diagrams}
\label{sec:pers-barc-diagr}

Throughout there is an underlying ideal of what persistent homology should be computing, which the field as a whole agrees on: given a filtered (and parametrized) sequence of topological spaces $\X_*$, the persistent homology $H_j^{a\to b}(\X_*)$ is the image of the induced map $H_j(\X_a)\to H_j(\X_b)$. In nice enough cases the collection of all such homologies has a nice algebraic description as some sort of collection of intervals, and these intervals with their start- and end-parameters can be used to produce diagrams that allow reasoning about the original spaces.

There are two main such diagrams in use -- both can be seen in Figure~\ref{fig:functionpersistence}. One view is the \dfn{persistence barcode} -- the sequence of interval is drawn, stacked on top of each other. Such a barcode can be seen in the middle of Figure~\ref{fig:functionpersistence}. The rank of any particular $H_j^{a\to b}(\X_*)$ is the number of intervals in the barcode that entirely covers the interval $(a,b)$. 

The other diagram in use is the \dfn{persistence diagram}: the start- and end-points of an interval in the interval decomposition of the persistent homology are taken to be $x$- and $y$-coordinates of points in the upper half of the first quadrant of the plane. An example can be seen to the right of Figure~\ref{fig:functionpersistence}. The number of points contained in the quadrant delimited by the horizontal line at height $a$ and the vertical line at width $b$ determines the rank of $H_j^{a,b}(\X_*)$.

Either of these cases is a visualization of the underlying data of a \emph{barcode}, which we can define as \textcite{cohen-steiner_stability_2007} as a multiset in $\eR^2$. The barcode is usually taken to include the uncountably many points along the diagonal of $\eR^2$ as part of the barcode.

Several metrics have been proposed for the space of all such barcodes or diagrams -- most of them have definitions more easily handled by working with the persistence diagram definition. We shall meet a few in this paper. For their definitions we shall assume that $X$ and $Y$ are two barcodes. We write $\Bij(X,Y)$ for the collection of all bijections between $X$ and $Y$. For up to countably many non-diagonal barcode elements, such bijections may pair each non-diagonal element with some possibly diagonal element, and pair all diagonal elements with infinitesimally close diagonal elements. 

The first two definitions here are taken from \cite{cohen-steiner_stability_2007}. The definition of Wasserstein distance is from \cite{cohen-steiner_lipschitz_2010}.

\begin{definition}
  The \dfn{bottleneck distance} $d_B(X,Y)$ is defined as
  \[
  d_B(X,Y) = \inf_{\gamma\in\Bij(X,Y)}\sup_{x\in X}\|x-\gamma(x)\|_\infty\quad.
  \]
\end{definition}

\begin{definition}
  The \dfn{Hausdorff distance} $d_H(X,Y)$ is defined on multisets $X,Y$ in $\R^2$ by
  \[
  d_H(X,Y) = \max\{\sup_{x\in X}\inf_{y\in Y}\|x-y\|_\infty, \sup_{y\in Y}\inf_{x\in X}\|y-x\|_\infty\}\quad.
  \]
\end{definition}

\begin{definition}
  The \dfn{Wasserstein distance} $d_W^p(X,Y)$ is defined as 
  \[
  d_W^p(X,Y) = \left(\inf_{\gamma}\sum_x\|x-\gamma(x)\|_\infty^p\right)^{1/p}
  \]
  For diagrams from persistent homology, where the points come in different dimension, the total Wasserstein distance sums the infima for each dimension separately before computing the $p$-th root.
\end{definition}

\section{Functions on a manifold}
\label{sec:functions-manifold}

The study of persistent homology originates from \textcite{edelsbrunner_topological_2000}, who first define the term and provide an algorithm for the computation of persistent homology. Taking their inspiration from $\alpha$-shapes, the authors assume that a filtered simplicial complex is provided as input, and produce a description of its persistent homology. In a slightly later paper, \textcite{edelsbrunner_hierarchical_2001} demonstrate that persistent homology can be applied to morse complexes from piecewise linear functions on a manifold -- the filtered simplicial complex required is given by combining the morse complex cells with the function values at the critical points witnessing each cell.

From this point and onwards, one strongly present culture in persistent homology remains focused on the role of a function defined on a manifold as the input data for the method. This viewpoint has proven remarkably fruitful in the study of \emph{stability}, and provides the best tools we currently have for justifying topological inferences with persistent homology.

\begin{figure}
  \centering
  \includegraphics[width=0.5\textwidth]{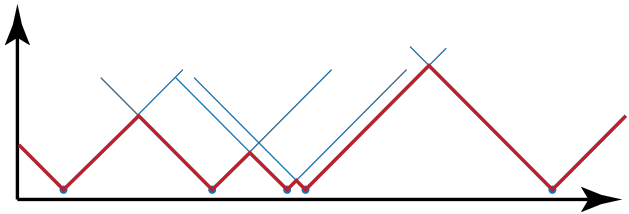}
  \caption{The distance to a set of points, defined for any point as the infimum of individual distances to points in the point cloud, produces a function for use in the functional approach to persistent homology. The points at the bottom of the valleys in the graph are the points of a 1-dimensional point cloud; and the lightly drawn cones emanating from each point correspond to the distance function from that point itself. The lower envelope of these distances forms the distance to the entire set, thus the function for encoding \v Cech complex persistent homology as a functional persistent homology.}
  \label{fig:pointcloudfunctional}
\end{figure}

It is worth noting that a point cloud topology point of view fits in this framework: as is illustrated in Figure~\ref{fig:pointcloudfunctional}, the distance to a discrete set of points produces a real-valued function on the ambient space of the points, with a persistent homology corresponding closely to the \v Cech complex homology of the point cloud itself.

\subsection{A functional view of persistent homology}
\label{sec:funct-view-pers}

\begin{figure}
  \includegraphics[page=1,width=0.7\textwidth]{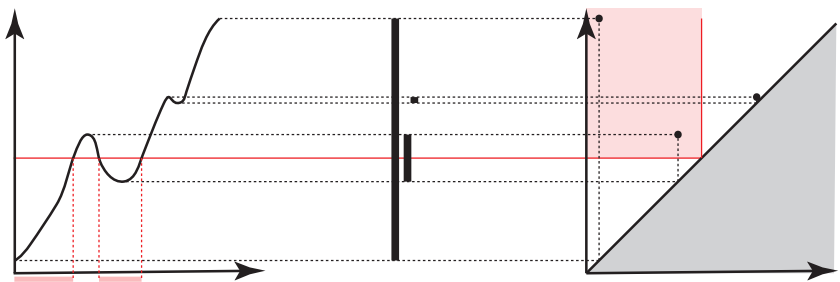}
  \caption{Persistence of $H_0$ of sublevel sets of a function $\mathbb R\to\mathbb R$. In black, we see the three components that appear at different times show up -- in the middle in a persistence barcode, and to the right as the three points in a persistence diagram. In red, we indicate a particular choice of height $\varepsilon$, at which the sublevel set has two components -- drawn below the graph to the left. These two components can be read off in both persistence visualizations -- through the two intersected bars in the middle, and through the two points contained in the shaded red region to the right.}
  \label{fig:functionpersistence}
\end{figure}

With this viewpoint, the fundamental given datum is a geometric object $\mathbb X$ and a tame function $f: \mathbb X\to\mathbb R$. In order to study the behavior of sublevel sets of $f$, persistent homology is used to measure the filtration of $\mathbb X$ given by $\mathbb X_\varepsilon = f^{-1}((-\infty, \varepsilon])$. 

A function $f:\mathbb X\to\mathbb R$ is called \emph{tame} if it is continuous, all sublevel sets have homology groups of finite rank, and there are finitely many critical values where the homology groups change.

This viewpoint, and the reasons for some of the choices made in creating algorithms are at their most apparent when considering the 1-dimensional case, where $\mathbb X=\mathbb R$, and we consider sublevel sets of some function $\mathbb R\to\mathbb R$.

Consider Figure~\ref{fig:functionpersistence}. Critical points of the function correspond to points where the sublevel set topology changes -- at minima, a new component is born, and at maxima, two components merge. To reflect these correspondences, we pair up critical points, choosing to pair a maximum with the latest relevant minimum, to reflect that the newer connected component merges in with the older one. The red line gives an example of a particular choice of height; the sublevel sets are split into two components, a fact reflected in the two bars intersected by the red line in the barcode -- the number of bars at any given parameter value reflects the corresponding Betti number at that stage.

We write $\Dgm_p(f)$ for the collection $\{(b, d)\}$ of start and endpoints of the barcode corresponding to the $p$th persistent Betti number $\beta_p$ of $f:\mathbb X\to\mathbb R$. 

We notice that the filtered simplicial complexes described in Section~\ref{sec:foundations-use} are the natural representations when the function studied is the distance to the sampled point cloud.

\subsection{Filtered complexes}
\label{sec:filtered-complexes}

The original persistence algorithm was formulated in terms of filtered complexes, and the functional view is fast to generate a filtered complex from the function under study. The key method to do this is described in \textcite{edelsbrunner_hierarchical_2001}: in a Morse theory approach, cells of a cellular complex correspond to critical values of the function, and depending on the index of the critical point, we can read off the dimensionality of the cell.

\begin{figure}
  \centering
  \includegraphics[width=0.7\textwidth]{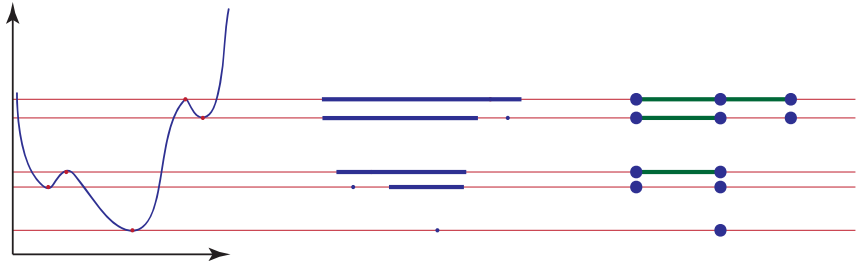}
  \caption{Going from a function on a manifold to a filtered sequence of spaces. Vertices of the Morse complex are given by the local minima, and each local maximum witnesses an edge connecting two neighbouring minima. To the left, we see the function with critical points marked, in the middle the sublevel sets at these points, and to the right, the corresponding filtered Morse complex. The filtered and parametrized structure of both spaces and complexes is clearly visible.}
  \label{fig:functiontofiltration}
\end{figure}

The Morse theoretic viewpoint gives a translation dictionary between critical points and cells in all dimensions, even where the example given in Figure~\ref{fig:functiontofiltration} is working in just one dimension. The fundamental feature to pay attention to is the \emph{index} of a critical point -- the number of negative signs in the appropriate quadratic form formulation of the Hessian at the critical point -- the higher the index, the higher the dimension of the cell corresponding to that critical point and introduced at the parameter of its function value in a sublevel set filtration.

\subsection{The stability meta-theorem}
\label{sec:stability}

There are results that depend crucially on a functional viewpoint in order to even articulate the question much less reach an answer. Most important of these is the issue of \emph{stability}. An introductory description would be that stability produces a continuity guarantee for the process that goes from a function to a barcode or persistence diagram descriptor of its persistent homology. If we can bound the change of a function, the resulting topological description should have bounded variation.

Stability theorems have the following general shape 
\begin{theorem}[Stability meta-theorem]
  For a \emph{nice enough space} $\X$ and \emph{nice enough functions} $f,g:\X\to\R$, a \emph{nice enough norm} of the difference $f-g$ is an upper bound to the distance between the barcodes of $f$ and $g$ in some \emph{nice enough metric}.
\end{theorem}

Most of the energy going into the study of stability has been improving these concepts of \emph{nice enough}, with significant and useful results. The development has relied at several stages on developing appropriate algebraic foundations to enable better theorem statements and more generous stability results.

\subsection{Vector spaces with ordered bases}
\label{sec:vector-spaces-with}

The first algebraic foundation in use was to consider the result and the intermediate computational stages of the persistence algorithm to be a vector space with a particular and ordered basis chosen. This viewpoint is implicit in \textcite{edelsbrunner_topological_2000}, where it generates the first algorithm for computing persistent homology.

\subsubsection{Persistence diagrams and stability}
\label{sec:pers-diagr-stab}

The work by \textcite{cohen-steiner_stability_2007} proves the first stability theorem for persistent homology: for a collection of persistent homology groups (referring to \cite{edelsbrunner_topological_2000} for their definition), the authors prove:
\begin{theorem}\label{thm:stability-1}
  Let $\X$ be a triangulable space with continuous tame functions $f, g:\X\to\R$. Then the persistence diagrams satisfy $d_B(\Dgm_p(f),\Dgm_p(g))\leq\|f-g\|_\infty$.
\end{theorem}

Here, tame is defined to mean that the function $f$ has a finite number of homological critical values and that the homology groups of sublevel sets are all finite-dimensional. This theorem was first proven, restricted to $p=0$, by \textcite{damico_optimal_2003}, using the language of \emph{size theory}.

The results from \textcite{cohen-steiner_stability_2007} have been generalized in several steps since its publication. Many if not most of these generalizations include a variation in the algebraic foundations to enable their greater power of proof.

\subsection{Diagrams over $(\R,\leq)$}
\label{sec:diagrams-over-r}

The first paper generalizing the results in \cite{cohen-steiner_stability_2007} was published by \textcite{chazal_proximity_2009}. The paper defines a persistence module $\mcF$ to be a diagram in the category $\Vect_\kk$ of the shape of the total order $(\R,\leq)$. In other words, $\mcF$ assigns a vector space $\mcF(x)$ to each $x\in\R$, and a linear map $\mcF(x\leq y)$ to each order relation $x\leq y$, making $\mcF$ a functor $(\R,\leq)\to\Vect_\kk$. We shall refer to these persistence modules as $(\R,\leq)$-modules and to the maps $\mcF(x\leq y)$ as translation maps.. The authors define a new tameness notion, and are able to prove an extended stability theorem.

\begin{definition}
  A $(\R,\leq)$-module $\mcF$ is \dfn{$\delta$-tame} if for any $\alpha<\alpha+\delta<\beta$ the rank of $\mcF(\alpha\leq\beta)$ is finite.
\end{definition}

\begin{definition}
 A function $f:\X\to\R$ is said to be $\delta$-tame if the $(\R,\leq)$-module of the homologies of the sublevel set filtration of $\X$ generated by $f$ is $\delta$-tame.
\end{definition}

The authors also define weak and strong interleaving -- concepts that will re-surface repeatedly in this direction of study. The articulation of the original definitions will be easier if we write $\X_f(\alpha)$ for the set $f^{-1}((-\infty),\alpha])$. 

\begin{definition}
  Two functions $f,g:\X\to\R$ are \dfn{weakly $\varepsilon$-interleaved} for $\varepsilon>0$ if there is some $a\in\R$ such that 
  \[
  \X_f(a+2n\varepsilon)\subseteq
  \X_g(a+(2n+1)\varepsilon)\subseteq
  \X_f(a+2(n+1)\varepsilon)
  \]
  for all $n\in\Z$. 
\end{definition}

\begin{definition}
  Two functions $f,g:\X\to\R$ are \dfn{strongly $\varepsilon$-interleaved} for $\varepsilon>0$ if for all $a\in\R$,  
  \[
  \X_f(a)\subseteq\X_g(a+\varepsilon)\subseteq\X_f(a+2\varepsilon)\quad.
  \] 
\end{definition}

Equivalently, the exact same formulas as for weak interleaving can be used, with the provision that they hold for all $a\in\R$, not just for some specific choice.

The definitions extend directly to generic $(\R,\leq)$-modules by the following definition:
\begin{definition}
  Let $\mcF$ and $\mcG$ be two $(\R,\leq)$-modules. Consider the diagram

  \begin{tikzpicture}
    \matrix (m) [matrix of math nodes, row sep=2em, column sep=4em]
    {
      \cdots & \mcF(a+2n\varepsilon) & 
      \mcF(a+(2n+1)\varepsilon) & \mcF(a+2(n+1)\varepsilon) & \cdots \\
      \cdots & \mcG(a+2n\varepsilon) & 
      \mcG(a+(2n+1)\varepsilon) & \mcG(a+2(n+1)\varepsilon) & \cdots \\
    };
    \foreach \i in {1,2} {
      \foreach \j/\jj in {1/2,2/3,3/4,4/5} {
        \draw [->] (m-\i-\j) -- (m-\i-\jj);
      }
    }
    \draw [->] (m-2-1) -- (m-1-2);
    \draw [->] (m-1-2) -- (m-2-3);
    \draw [->] (m-2-3) -- (m-1-4);
    \draw [->] (m-1-4) -- (m-2-5);
  \end{tikzpicture}

  The modules $\mcF$ and $\mcG$ are \dfn{weakly $\varepsilon$-interleaved} if there are linear maps for the diagonal arrows that make all these diagrams commute for some fixed $a\in\R$ and all $n\in\Z$.

  The modules $\mcF$ and $\mcG$ are \dfn{strongly $\varepsilon$-interleaved} if there are linear maps for the diagonal arrows that make all these diagrams commute for all $a\in\R$ and all $n\in\Z$.  
\end{definition}

With these definitions in place, we are able to state the most important results from \cite{chazal_proximity_2009}.

\begin{proposition}
  Let $f,g:\X\to\R$ be two real-valued functions on a topological space.
  \begin{enumerate}
  \item If $f,g$ are strongly $\varepsilon$-interleaved, then they are weakly $\varepsilon$-interleaved.
  \item If $f,g$ are weakly $\varepsilon$-interleaved, then they are strongly $3\varepsilon$-interleaved.
  \item $f,g$ are strongly $\varepsilon$-interleaved if and only if $\|f-g\|_\infty\leq\varepsilon$.
  \end{enumerate}
\end{proposition}

The last of these three parts is crucial for this approach to stability -- it means that stability results can be translated into how persistence diagrams of interleaved modules behave rather than how persistence diagrams depend on functional properties.

For a persistence diagram $D$, we can define the $\delta$-persistence diagram $D_\delta$ by removing any points within $\delta$ of the diagonal, i.e. any point $(b,d)$ with $0<d-b<\delta$. Instead of augmenting the diagram with the diagonal $\Delta=\{(x,x) : x\in\R\}$, we augment $D_\delta$ with the diagonal $\{(x,x+\delta): x\in\R\})$. 

Write $\mcI_a^b$ for the $(\R,\leq)$-module that has $\mcI_a^b(x)=\kk$ if $x\in[a,b]$ and $\mcI_a^b(x)=0$ otherwise, and such that all $\mcI(x\leq y)$ are the zero map unless $a\leq x\leq y\leq b$, in which case $\mcI_a^b(x\leq y)$ is the identity map. We call $\mcI_a^b$ the \dfn{interval module} for the interval $[a,b]$. Similar definitions can be produced for $(a,b]$, $[a,b)$ and $(a,b)$. Some $(\R,\leq)$-modules decompose into a direct sum of interval modules.

For an $(\R,\leq)$-module $\mcF$ that does decompose into an interval module, we write $\Dgm(\mcF)$ for the multiset of endpoints of intervals in a decomposition in such modules. Since the interval modules are indecomposables in the category of $(\R,\leq)$-modules, it follows from the Krull-Schmidt-Azumaya-theorem that the decomposition, and therefore this diagram, is unique if it exists.

\begin{theorem}
  Suppose $\mcF$ and $\mcG$ are weakly $\varepsilon$-interleaved $(\R,\leq)$-modules that both are $\delta$-tame for some $\delta\geq0$. Then
  \[
  d_B(\Dgm(\mcF)_\delta,\Dgm(\mcG)_\delta) \leq 3\varepsilon
  \]
\end{theorem}

\begin{theorem}
  Suppose $\mcF$ and $\mcG$ are strongly $\varepsilon$-interleaved $(\R,\leq)$-modules that both are $\delta$-tame for some $\delta\geq0$. Then
  \[
  d_H(\Dgm(\mcF)_\delta,\Dgm(\mcG)_\delta) \leq 3\varepsilon
  \]
\end{theorem}

In particular it follows that 
\begin{theorem}\label{thm:stability-2}
  If $f,g$ are two $\delta$-tame functions such that $\|g-f\|_\infty\leq\varepsilon$ for $\delta\geq0$ and $\varepsilon>0$, then for any $p$,
  \[
  d_B(\Dgm_p(f)_\delta,\Dgm_p(g)_\delta) < \varepsilon
  \]
\end{theorem}

It is worth noticing that for Theorem~\ref{thm:stability-2}, the assumptions on triangulability for $\X$ and on continuity for $f, g$ from Theorem~\ref{thm:stability-1} have been removed.

\begin{theorem}
  Suppose $L$ is a finite point cloud in some metric space. There are $(\R,\leq)$-modules
  \begin{align*}
  H_p\vC(L)(a) &= H_p(\vC_{2^a}(L)) \\
  H_p\VR(L)(a) &= H_p(\VR_{2^a}(L)) \\
  H_p W(L,W)(a) &= H_p(W_{2^a}(L,W))
\end{align*}
with all the translation maps induced from the inclusion maps. 
  
  Then
  \[
  d_B(\Dgm(H_p\vC(L)),\Dgm(H_p\VR(L))) \leq 1\quad.
  \]

  If the points of $L$ are densely sampled from a compact set 
  $L\subseteq W\subseteq\X\subset\R^d$, with sampling conditions stated in \cite[Theorem 3.7]{chazal_towards_2008},
  \[
  d_B(\Dgm(H_p\vC(L)),\Dgm(H_p W(L,W))) \leq 3\quad.
  \]
\end{theorem}

\subsection{Directly to diagrams}
\label{sec:directly-diagrams}

Some of the stability results skip the intermediate step of persistence modules altogether, and argue entirely in terms of the persistence diagram and its behaviour. As far as we have been able to tell, this operates with an underlying assumption of using  $(\R,\leq)$-modules as an algebraic framework, but some papers never articulate this choice concretely.

\textcite{cohen-steiner_lipschitz_2010} prove
\begin{theorem}
  Let $\X$ be a triangulable compact metric space implying bounded degree $k$ total persistence, for $k\leq 1$, and let $f,g:\X\to\R$ be two tame Lipschitz functions. Then
  \[
  d_W^p(f,g) \leq C^{\frac{1}{n}}\cdot\|f-g\|_\infty^{1-\frac{k}{p}}
  \]
  for all $p\geq k$, where $C=C_\X \max\{\Lip(f)^k,\Lip(g)^k\}$.
\end{theorem}

Here, a metric space $\X$ implies bounded degree $k$ total persistence if there is some $C_\X$ depending only on $\X$ such that $\Pers_k(f)\leq C_\X$ where $D_f$ is the persistence diagram of a sublevel set filtration of $f:\X\to\R$ with Lipschitz constant $\Lip(f)\leq 1$ and $\Pers_k(f)=\sum_{(b,d)\in D_f, b-d>t}(b-d)^k$. We call a Lipschitz function \dfn{tame} if the homologies of the sublevel sets come with finitely many changes and each homology group has finite rank.

With the same notation, there is a stability theorem for the total persistence moments $\Pers_p(f)$ too. We write $\Amp(f)=\max_{x\in\X}f(x)-\min_{y\in\X}f(y)$.
\begin{theorem}
  Let $\X$ be a triangulable, compact metric space that implies bounded degree-$k$ total persistence for $k\geq 0$, and let $f,g:\X\to\R$ be two tame Lipschitz functions. Then
  \[
  |\Pers_p(f) - \Pers_p(g)| \leq 4pw^{p-1-k}C\cdot\|f-g\|_\infty
  \]
  for every $p\geq k + 1$, where $C=C_\X \max\{\Lip(f)^k,\Lip(g)^k\}$ and $w$ is bounded from above by $\max\{\Amp(f ), \Amp(g)\}$.  
\end{theorem}

\subsection{Measures on the real line}
\label{sec:measures-real-line}

The approach in Section~\ref{sec:diagrams-over-r} was elaborated by~\textcite{chazal_structure_2012}. In that paper, the authors deal primarily with the fundamental question of which conditions on a total order module allow for a persistence diagrams decomposition to even exist. For cases where these decompositions do exist, they are able to prove stability theorems; and in order to establish this existence, they develop a fruitful notation and viewpoint.

Their approach continues with the emphasis on the behaviour of \emph{persistence diagrams} that we saw in Section~\ref{sec:directly-diagrams}. There are theorems in their work that relates the work to the behaviour of specific parametrized filtrations on concrete spaces, but most of the work considers persistence diagrams of abstracted and decomposable $(\R,\leq)$-modules with specific tameness conditions directly.

\subsubsection{Persistence measures}
\label{sec:persistence-measures}

At the core of this approach is the recognition that multisets of points in the plane correspond precisely to locally finite integer-valued measures on the plane, that can then be considered to be counting the points as point masses. To elaborate, the authors consider four types of persistence intervals: $[a,b], [a,b), (a,b], (a,b)$. To acquire a coherent notation for these, they introduce point decorations -- $a^+$ can be thought of $a+\varepsilon$ for some infinitesimal $\varepsilon$, so that an interval starting at $a^+$ is open at that end, and an interval ending in $a^+$ is closed at that end. Similarily, $b^-$ can be thought of as $b-\varepsilon$ for an infinitesimal $\varepsilon$, so that an interval starting in $b^-$ is closed, and an interval ending in $b^-$ is open. Following~\cite{chazal_structure_2012}, we write $a^*$ when we do not have any information about the decoration of $a$. 

Viewing an interval $(a^*, b^*)$ as a point in a persistence diagram -- viewed as a multiset in the plane -- the point is some $(a,b)$ decorated with a flag pointing in one of the quadrant directions: $++, +-, -+,$ or $--$. Translation from a persistence diagram to a \emph{persistence measure} now follows easily: for a rectangle $[a,b]\times[c,d]$ in the plane (with undecorated endpoints), a point is counted by the measure if the point flag points into the interior of the rectangle. With this definition, additivity for the measure function can be proven, and a whole slew of measure theoretic machinery can be used. 

Far from all possible $(\mathbb R,\leq)$-modules are decomposable into interval modules -- \textcite{de_silva_persistence_2012} gives as examples of decomposable total order modules the classes
\begin{enumerate}
\item Modules over finite orders (proven by~\textcite{gabriel_unzerlegbare_1972}).
\item Modules over $(\mathbb Z, \leq)$ of locally finite dimension (proven by~\textcite{webb_decomposition_1985}).
\item Modules over $(\mathbb R, \leq)$ of locally finite dimension (proven by~\textcite{crawley-boevey_decomposition_2012}).
\end{enumerate}
but also points out that~\textcite{webb_decomposition_1985} demonstrates a module $M$ over $(\mathbb R,\leq)$ where each $M(x\leq y)$ has finite rank, but still the module is not decomposable into intervals.

The measure approach, however, provides a decomposition into interval modules whenever one actually exists, and if no global decomposition exists, can identify regions of the plane where the persistence diagram is decomposable, and provide a decomposition over these regions.

For an arbitrary representation of $(\mathbb R, \leq)$, the authors are able to prove a measure $\mu$ that coincides with the point mass perspective: for a persistence module $M$ such that every $M(x\leq y)$ has finite rank $r_x^y$, the measure $\mu([a,b]\times[c,d])$ is defined to be $r_b^c-r_a^c-r_b^d+r_a^d$. This way, the measure $\mu$ is defined for arbitrary modules with finite rank translation maps. 

\subsubsection{Tameness conditions}
\label{sec:tameness-conditions}

Based on these decomposability regions, the authors are able to define a family of tameness conditions, with inclusions of classes of modules along the arrows:

\begin{center}
\begin{tikzpicture}[->,thick,node distance=2cm]
  \matrix (m) [matrix of nodes, row sep=2ex, column sep=4em] 
  {
    & & & v-tame & \\
    finite & locally finite & q-tame & & r-tame \\
    & & & h-tame & \\
  };
  \draw (m-2-1) -- (m-2-2);
  \draw (m-2-2) -- (m-2-3);
  \draw (m-2-3) -- (m-1-4);
  \draw (m-2-3) -- (m-3-4);
  \draw (m-1-4) -- (m-2-5);
  \draw (m-3-4) -- (m-2-5);
\end{tikzpicture}
\end{center}

\begin{figure}
  \centering
  \begin{tikzpicture}[scale=4]
    \fill [gray] (0,0) -- (1,1) -- (0,1) -- cycle;
    \fill [gray] (0,0) -- (-0.2,0) -- (-0.2,1.2) -- (0,1.2) -- cycle;
    \fill [gray] (-0.2,1) -- (-0.2,1.2) -- (1,1.2) -- (1,1) -- cycle;
    \filldraw [black, fill=white,opacity=0.7] (0.3,0.6) rectangle (0.5,0.8);
    \node at (0.4,0.7) {$R$};
    \filldraw [black, fill=white,opacity=0.7] (0.3,0.81) rectangle (0.5,1.2);
    \node at (0.4,0.9) {$V$};
    \filldraw [black, fill=white,opacity=0.7] (-0.2,0.6) rectangle (0.29,0.8);
    \node at (0.2,0.7) {$H$};
    \filldraw [black, fill=white,opacity=0.7] (-0.2,0.81) rectangle (0.29,1.2);
    \node at (0.2,0.9) {$Q$};
    \fill [white] (0,-0.1) -- (0,1.3) -- (-0.1,1.3) -- (-0.1,-0.1) -- cycle;
    \fill [white] (-0.3,1) -- (1.1,1) -- (1.1,1.1) -- (-0.3,1.1) -- cycle;
  \end{tikzpicture}
  \caption{Tameness-conditions schematically illustrated.}
  \label{fig:tameness}
\end{figure}
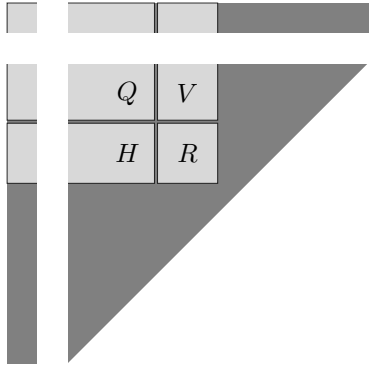

Here, a module $M$ is...
\begin{description}
\item[finite] if $M$ is a finite direct sum of interval modules.
\item[locally finite] if $M$ is a direct sum of interval modules, such that only finitely many span any given $t\in\mathbb R$.
\item[q-tame] if the measure corresponding to $M$ is finite over every quadrant not touching the diagonal.
\item[h-tame] if the measure corresponding to $M$ is finite over every horizontally infinite strip $H$ not touching the diagonal.
\item[v-tame] if the measure corresponding to $M$ is finite over every vertically infinite strip $V$ not touching the diagonal.
\item[r-tame] if the measure corresponding to $M$ is finite over every finite rectangle not touching the diagonal.
\end{description}
These four last cases are sketched out in Figure~\ref{fig:tameness}; the horizontal bars above and to the left correspond to interval modules that survive until $+\infty$ and interval modules that were born at $-\infty$ respectively. 

\subsubsection{Order module view of interleaving}
\label{sec:order-module-view}

The plane $\bar{\mathbb R}^2$ has a partial order given by $(p_1,q_1)\leq(p_2,q_2)$ if and only if $p_1\leq q_1$ and $p_2\leq q_2$. We can define the shifted diagonals $\Delta_x = \{(p,q) | q-p=2x\}$ as subsets of the plane; with order structure inherited from this order on the plane. These diagonals are isomorphic -- as posets -- to $(\mathbb R,\leq)$: by picking $t\mapsto(t-x,t+x)$, this isomorphism is canonical.

With this structure, the authors define (strong) $|y-x|$-interleaving of persistence modules $M, N$ as the existence of a persistence module $I$ over $\Delta_x\cup\Delta_y$ such that $I|_{\Delta_x}=M$ and $I|_{\Delta_y}=N$. The authors prove that if $M, N$ are $\delta$-interleaved, then there is a family of persistence modules $P_x$ such that $P_0=M$, $P_\delta=N$, and $P_x, P_y$ are $|y-x|$-interleaved for all $x,y\in[0,\delta]$. These can be fused together into a single persistence module over the diagonal strip $\{(p,q): 0\leq q-p\leq 2\delta\}$ with the above partial order structure.

From this interleaving relation, the authors define an interleaving distance: 
\[
d_i(M, N) = \inf\{\delta : M, N \text{ are $\delta$-interleaved} \}
\]

This distance is a pseudo-metric: the authors prove the triangle inequality, but give as an example the four interval modules for the intervals $(p^-,q^-),(p^+,q^-),(p^-,q^+),(p^+,q^+)$, that all have 0 interleaving distance, but are not in fact isomorphic.

\subsubsection{Stability}
\label{sec:stability-1}

Based on this machinery, the authors are able to prove a number of stability-related theorems, that all lead to the fundamental isometry theorem, occuring in \cite[Theorem 4.11]{chazal_structure_2012}, and also proven independently by~\textcite{lesnick_optimality_2011}:

\begin{theorem}[Isometry]
  Let $M, N$ be q-tame persistence modules. Then $d_i(M,N)=d_b(\Dgm(M),\Dgm(N))$.
\end{theorem}

The applications of the isometry come from identifying tameness conditions for classes of persistent homology modules:

\begin{theorem}[Theorem 2.23 of~\cite{chazal_structure_2012}]
  Let $X$ be a locally compact polyhedron, and $f:X\to\mathbb R$ a proper continuous function. Then the persistent homology of the sublevel set filtration of $(X, f)$ is h-tame, v-tame, and r-tame, but not q-tame.
\end{theorem}

Notice that the collection of tameness conditions that hold here mean that as long as we ignore any parts of the persistent homology of the sublevel set filtration that persists all the way from $-\infty$ to $\infty$, the remainder is tame enough for the isometry theorem, and therefore stable.

\begin{theorem}[Proposition 5.1 of~\cite{chazal_persistence_2012}]
  If $(X, d_X)$ is a precompact metric space, then the \v Cech and Vietoris-Rips persistent homology modules are q-tame.
\end{theorem}

It is also well-known in the community that strong finiteness conditions, and therefore also q-tameness, hold for the homologies of sublevel filtrations if 
\begin{itemize}
\item $X$ if a compact manifold and $f$ is a Morse function.
\item $X$ is a compact polyhedron and $f$ is piecewise linear.
\end{itemize}

From the isometry theorem also follows, by the view of interleaving as a persistence module, the classical stability theorem of~\textcite{cohen-steiner_stability_2007} as we have already discussed in Section~\ref{sec:pers-diagr-stab}, Theorem~\ref{thm:stability-1}.

\subsection{Categorification}
\label{sec:categorification}

Work by \textcite{bubenik_categorification_2012} studies the category of functors $(\R,\leq)\to\Vect_\kk$, and is able to prove that the category of persistence modules is abelian. 

They leverage this to prove a generous stability theorem: for \emph{arbitrary} (not necessarily continuous) functions $\X\to\R$ from a topological space, and any functor $H$ from topological spaces to a category of real-indexed diagrams in an abelian category $\mathcal D$, the \emph{interleaving distance between the diagrams} generated by applying $H$ to the sublevel set filtrations of the functions is bounded above by the \emph{$L_\infty$-distance} of the functions. 

Furthermore, they prove many of the categories that emerge naturally in persistent homology are abelian.

\newpage
\section{Filtered topological spaces}
\label{sec:filt-topol-spac}


The other culture present in the study of persistent homology focuses on the role of a \emph{filtered topological space} and derived algebraic objects as the fundamental notion. This viewpoint has sparked a wealth of algebraic abstractions and given rise to several different notions of the \emph{shape} of a persistent homology theory.

Connecting this viewpoint with the original study by \textcite{edelsbrunner_topological_2000}, and indeed with the entire viewpoint present in Section~\ref{sec:functions-manifold}, one may point out that for any function $f:\X\to\R$, the sublevel sets $f^{-1}((-\infty,x])$ form a filtration of $\X$. For tame enough -- finitely many topological critical points, finite rank homology for any sublevel set, and similar conditions -- functions, this filtration can be described by a finite filtration, or even a parametrization with finitely many different states.

Since homology is a functor, and inclusions are continuous maps, applying homology to a filtration produces a diagram of homology groups on the shape
\[
H_j\X_0 \to H_j\X_1 \to \dots \to H_j\X_n
\]
and by interpreting this diagram as a module in one of a number of different possible module categories, further generalizations are possible.

Commonly, the geometric filtrations in use in persistent homology really are parametrizations -- for any value $\varepsilon\in\R$, there is some resulting space $\X_\varepsilon$ -- that happen to generate filtrations: 
\[
\X_{(-\infty,\varepsilon]} = 
\mkern-18mu \bigcup_{\delta\in(-\infty,\varepsilon]} \mkern-18mu X_\delta
\qquad\qquad
\varepsilon<\varepsilon' \Rightarrow 
\X_{(-\infty,\varepsilon]} \subseteq \X_{(-\infty,\varepsilon']}
\]
For these cases, it is common to blur the lines between the definitions of filtered spaces and parametrized spaces.

\subsection{Vector space with ordered basis}
\label{sec:vector-space-with}

If a simplicial complex is filtered, then this induces a preorder on simplices of the simplicial complex -- any simplex precedes all simplices from later filtration stages. Any preorder can be specialized to a total order by picking arbitrarily some ordering of elements that do not already have an ordering setup by the preorder -- and this is certainly the case with the preorder from a filtration. This total order can even be picked to be compatible with the coface relation on simplices.

In particular this means that from a filtered simplicial complex, we can easily construct a chain complex with a totally ordered simplex basis. This was the basis of the original algorithm in \cite{edelsbrunner_topological_2000}: simplices are consumed from a totally ordered stream, and the change in topology resulting from the inclusion of any one simplex is reflected in a changing state, from which barcodes can be read off.

This setting has also informed extensions to the work by \textcite{edelsbrunner_topological_2000}: in a paper by \textcite{cohen-steiner_extending_2009}, a total ordering of the simplices in a simplicial complex $K$ is used to filter $K$ both by taking initial sequences $K_j=(\sigma_0,\dots,\sigma_j)$ and by taking terminal sequences $L_j=(\sigma_j,\dots,\sigma_N)$ of the simplices. With these building blocks, then, the original persistent homology sequence
\[
H_*(K_0) \to \dots \to H_*(K_N)
\]
can be extended by taking homology relative to terminal sequences to produce an \dfn{extended persistence} sequence
\[
H_*(K_0)\to\dots\to H_*(K_N)\to H_*(K_N,L_0)\to\dots\to H_*(K_N,L_N)
\]

This sequence, motivated by Poincar\'e and Lefschetz duality, carries a number of benefits over the original persistent homology theory. One of them is that no infinite barcodes occur -- any interval will have an endpoint, possibly among the relative homology groups.

Duality produces numerous symmetry relations in the persistence diagram for extended persistence. The paper \cite{cohen-steiner_extending_2009} describes a way to draw the persistence diagram for the extended case so that the symmetries emerge as mirror symmetries in the diagram -- with some adjustments for dimension shifts in the duality results.

For the case where the ordering of cells comes from \emph{piecewise linear functions} on a \emph{simplicial complex}, extended persistence has a stability theorem bounding the \emph{bottleneck distance} of diagrams by the \emph{$L_\infty$-norm of the difference between the corresponding functions}. See Section~\ref{sec:stability} for more details on stability.

Quite some research has gone into optimizing the persistence algorithm in various ways. Here, handling the sorted boundary matrix tends to be at the center of attention. Some notable results include:

\textcite{milosavljevic_zigzag_2009} demonstrate a method using matrix multiplication to compute persistent homology; proving the existence of a method that runs in matrix multiplication time. This gives the currently best worst-case complexity estimate for computing persistent homology as $O(n^\omega)$ for $n$ simplices, and $\omega$ the matrix multiplication exponent. At the point of writing this article, $\omega$ is practical for $\log_2 7\approx 2.807$ \cite{strassen1969gaussian}, and possible, with a large constant, at $2.3727$ \cite{coppersmith1990matrix,williams2012multiplying}.

\textcite{cohen-steiner_vines_2006} demonstrate how the change in the persistence diagram induced by the re-ordering of simplices in the filtration can be traced in linear time: the result is a \emph{vineyard}, tracing the change in the persistence diagram induced by a homotopy between functions inducing filtrations of a simplicial complex. This approach allows the proof of a combinatorial stability theorem, bounding \emph{bottleneck distance} by the \emph{$L_\infty$-distance} between \emph{simplicial approximations of continuous functions} on a \emph{simplicial complex}.

\textcite{cohen-steiner_persistent_2009} study a functional persistence situation where $f:\X\to\R$, $g:\Y\to\R$, and $\Y\subseteq\X$. In this case, there is an induced map from the homology of the sublevel sets of $g$ to the homology of the sublevel sets of $f$, and the authors give algorithms for computing kernels, images, and cokernels of this induced map. Their algorithms fundamentally work with adjusting sorted matrices of simplices.

\textcite{chen_output-sensitive_2011} work with Monte Carlo algorithms for estimating the rank of a matrix to compute persistence barcodes from a sorted boundary matrix. This approach speeds up the computation of barcodes.

\textcite{chen_persistent_2011} notice that since the persistence algorithm induces a pairing between columns in the boundary matrix, pairing up a completely emptied out column with one that corresponds to the last simplex to bound the cycle, a re-ordering of the computation can eliminate many matrix operations. Their approach wcannot be modified to output concrete representative cycles, but improves asymptotic bounds for the problem of computing a barcode.

\subsection{Graded modules over $\kk[t]$}
\label{sec:graded-modules-over}

The first significant advance in the choice of underlying algebraic structure for persistence modules came from \textcite{cz2005}. They observe that a diagram of vector spaces
\[
V_0 \to V_1 \to \dots 
\]
can be modelled as a graded module over the polynomial ring $\kk[t]$. The module $V_*$ is taken to have $V_d$ in degree $d$, and the action of multiplying by $t$ corresponds to the linear map $V_d\to V_{d+1}$. Homology groups with field coefficients are vector spaces, and the inclusions of complexes induce linear maps between these vector spaces. Thus, this construction translates a persistent homology diagram to a graded module over $\kk[t]$.

At this stage, \textcite{cz2005} observe that the existence of a barcode decomposition follows directly from the fact that $\kk[t]$ is a principal ideal domain, and therefore any module $V_*$ decomposes into a direct sum of cyclic modules. These come in two versions: torsion modules isomorphic to $\kk[t]/(t^d)$ for some natural number $d$, and free modules isomorphic to $\kk[t]$. These two classes can be directly translated into free and finite intervals $[a,a+d)$ or $[a,\infty)$.

The work in \cite{cz2005} also demonstrates that the persistence algorithm described by \textcite{edelsbrunner_topological_2000} works with the same result for arbitrary field coefficients where the original description required coefficients in the field $\Z/2\Z$.

This work has been extensively cited -- to the point where the papers \cite{edelsbrunner_topological_2000,cz2005} are the standard reference citations for the persistence algorithm, and a number of extensions to the results have been provided, as well as numerous applications to the extension of expressive power the change of fields produces.

\subsubsection{Results relying on non-binary fields}
\label{sec:results-relying-non}

The most obvious direct usefulness of the graded polynomial ring module approach has been in cases where the dependency of homology on the characteristic of the coefficients matters. This was the case in work by \textcite{imagespace2008}.

 A study by \textcite{lee_nonlinear_2003} investigates the statistics of $3\times3$ pixel patches from naturally occuring images. They find, inter alia, a high density circle in the first few PCA coordinates. This circle, they notice, corresponds closely to linear gradient directions within the dataset.

\textcite{imagespace2008} pick up the same dataset, and study it using persistent homology. They are able to recover two additional, secondary, high-density circle shapes within the dataset. These three circles combine to form a high-density 2-dimensional surface, which after computing persistent homology over both $\Z/2\Z$ and over $\Z/3\Z$ could be identified as the Klein bottle.

The ability to compute persistent homology with coefficients in $\Z/3\Z$ was crucial for this approach, and algorithmically dependent on the graded module over $\kk[t]$ approach to persistent homology.

\subsubsection{Multi-dimensional persistence}
\label{sec:multi-dimens-pers}

With inspiration from the several relevant parameters affecting the analysis in \cite{imagespace2008}, \textcite{carlsson_theory_2009} constructed \dfn{multidimensional persistence}. The underlying observation is that just as graded modules over $\kk[t]$ model singly parametrized topological spaces, adding more parameters corresponds to adding more variables to the polynomial ring. Hence, a $d$-dimensional parametrization can be modeled in a persistence way by working with graded modules over $\kk[t_1,\dots,t_d]$.

The multi-dimensional theory has problems -- chief among which is the lack of as useful a decomposition into a small and easy to describe class of indecomposables. The category of graded modules over $\kk[t_1,\dots,t_d]$ has no complete discrete invariant, but \textcite{carlsson_theory_2009} propose a discrete invariant -- the rank invariant -- turning out to be incomplete but useful.

The theory has been further studied since:

\textcite{carlsson_computing_2009} introduce Gr\"obner basis methods for computing multidimensional persistent homology, demonstrating that for \emph{one-critical multifiltrations}, the rank invariant can be computed in polynomial time. The translation process they use to recast the problem to a Gr\"obner basis computation has potential exponential blowup behaviours for the general case. \textcite{patriarca_presentation_2012} demonstrate that by avoiding the mapping telescope and using more refined Gr\"obner basis approaches, the computation can be bounded to polynomial time in general.

The multidimensional approach has received a lot of attention from the size function community. \cite{biasotti_multidimensional_2008,cerri_multidimensional_2009,cerri_betti_2010,cagliari_finiteness_2011} treat multidimensional persistent homology in a size function framework as important tools for image analysis. 

Questions of stability for multi-dimensional persistence modules have been studied, both in the size function community (\cite{cerri_multidimensional_2009,cerri_betti_2010}) and in the context of persistent homology by \textcite{lesnick_optimality_2011}.

\subsubsection{Cohomology and duality}
\label{sec:cohomology-duality}

Persistent cohomology was mentioned by \textcite{cohen-steiner_extending_2009}, who immediately use Lefschetz duality to transform it into relative homology. Independently, \textcite{de_silva_persistent_2009}, later extended by \textcite{morozov_persistent_2011}, produce an algorithm for computing persistent cohomology and observe connections to computing intrinsic circle-valued coordinate functions from point cloud datasets.

This work inspired a paper by \textcite{de_silva_dualities_2011} in which two duality functors -- $M_*\mapsto\hom_\kk(M_*,\kk)$ and $M_*\mapsto\hom_{\kk[t]}(M_*,\kk[t])$ on graded $\kk[t]$-modules are studied, and how these functors affect both the persistence algorithm itself, the ordering of basis elements in a sorted vector space approach, and how the barcodes are modified. These two functors allow the transport of information between relative and absolute versions of persistent homology and cohomology. 

\subsubsection{Algebraic adaptation of topological constructions}
\label{sec:algebr-adapt-topol}

In ongoing work, \textcite{lsvj2011parallelpersistence} work out algorithms and approaches for using spectral sequences of graded $\kk[t]$-modules to parallelize the computation of persistent homology. The approach fundamentally relies on the algebra of graded $\kk[t]$-modules as a proxy for persistent homology.

From this work, \textcite{skraba_algebraic_2013} work out more detailed algorithmics for graded $\kk[t]$-modules and are able to generalize the results from \textcite{cohen-steiner_persistent_2009} to allow computation of images, kernels, and cokernels of a wider range of maps in persistent homology.

\subsection{Modules over a quiver algebra}
\label{sec:modules-over-quiver}

Another algebraic model that describes the persistent homology diagrams of vector spaces is given by quiver algebras. A persistent homology diagram of the shape
\[
H_k(\X_0) \to H_k(\X_1) \to
H_k(\X_2) \to H_k(\X_3) \to
H_k(\X_4) \to H_k(\X_5)
\]
can be considered as a module over the path algebra $\kk Q$ for $Q$ the quiver
\[
\bullet \to \bullet \to 
\bullet \to \bullet \to 
\bullet \to \bullet
\]

A theorem by \textcite{gabriel_unzerlegbare_1972} states:

\begin{theorem}[Gabriel's theorem]
  A quiver $K$ has finitely many isomorphism classes of irreducible finite dimensional $\kk$-linear representations if and only if $K$ is a disjoint union of finitely many quivers of the classes $A_e$, $D_m$, or $E_n$ for $e\geq 1$, $m\geq 4$, $6\leq n\leq 8$. 
\end{theorem}

In particular, Gabriel goes on to prove that the exact isomorphism classes that show up for the quivers of type $A_e$ -- quivers of linear sequences of arrows, possibly alternating in direction -- are the interval modules. These have some connected interval along the linear sequence where one-dimensional vector spaces are connected by identity maps -- and outside this interval, all maps are and all vector spaces are zeros.

For the case of ``classical'' persistent homology, this recovers the barcode description for the case of a filtered finite simplicial complex: the persistent homology decomposes into a direct sum of irreducibles, and these irreducibles all are these interval modules. To describe each interval module, it is enough to state its start and end index, which is the exact data that a barcode conveys.

This approach has given rise to two generalization directions in particular.

\subsubsection{Zigzag persistence}
\label{sec:zig-zag-persistence}

\textcite{carlsson_zigzag_2008} pointed out that Gabriel's theorem has concrete consequences for topological data analysis. In particular, the non-dependency on arrow direction for a quiver to qualify as having type $A_e$ means that we can consider quivers where arrows alternate direction, either occasionally or consistently.

This paper introduces the fundamental idea, provides matrix algorithms for computing zigzag persistence, and provides the \emph{diamond principle}, relating how local changes along the zigzag reflect in changes to the persistence diagram. The paper also suggests several applications where the zigzag naturally arises:
\begin{description}
\item[Balancing different parameters] In the study by \textcite{imagespace2008}, the $p\%$ densest points as computed with a parametrized density estimator were used to determine the topology of the dataset. For studies like this one, it is worth while to try to work with all possible values of the parameter determining the density estimator at once -- to replicate the success persistent homology has in sweeping over entire ranges for a parameter. 

Writing $\X_r^p$ for the densest $p\%$ of the point cloud $\X$ as measured using the parameter $r$. Varying $r$ along $r_1<r_2<\dots<r_N$, there is a zigzag
  \[
  \begin{tikzpicture} [thick, ->]
    \matrix (m) [matrix of math nodes,row sep=2ex] 
    {
      & X_{r_1}^p \cup X_{r_2}^p &
      & X_{r_2}^p \cup X_{r_3}^p &
      & X_{r_3}^p \cup X_{r_4}^p &
      & \dots &
      & X_{r_{N-1}}^p \cup X_{r_N}^p  \\
      X_{r_1}^p &&
      X_{r_2}^p &&
      X_{r_3}^p &&
      X_{r_4}^p &&
      X_{r_{N-1}}^p &&
      X_{r_N}^p \\
    };
    \draw (m-2-1) -- (m-1-2);
    \draw (m-2-3) -- (m-1-2);
    \draw (m-2-3) -- (m-1-4);
    \draw (m-2-5) -- (m-1-4);
    \draw (m-2-5) -- (m-1-6);
    \draw (m-2-7) -- (m-1-6);
    \draw (m-2-7) -- (m-1-8);
    \draw (m-2-9) -- (m-1-8);
    \draw (m-2-9) -- (m-1-10);
    \draw (m-2-11) -- (m-1-10);
  \end{tikzpicture}
  \]
  For each point cloud in this sequence, compute a geometric complex, and compute its homology -- the resulting diagram is a zigzag diagram, and its decomposition into barcodes carries information about the variation of $r$ in a way directly analogous to how persistent homology itself measures homological features over varying values for a parametrization.
\item[Topological bootstrapping] Similar to bootstrapping in statistics, one may want to take a sequence of small samples $\X_i$ from a large dataset $\X$ and estimate the topology of each $\X_i$ individually. Doing this, disambiguation between local features of each $\X_i$ and global features of $\X$ is not entirely transparent.

Here, the \dfn{union zigzag} provides a method for persisting features across several samples:
  \[
  \begin{tikzpicture} [thick, ->]
    \matrix (m) [matrix of math nodes,row sep=2ex] 
    {
      & X_1 \cup X_2 &
      & X_2 \cup X_3 &
      & X_3 \cup X_4 &
      & \dots &
      & X_{N-1} \cup X_N  \\
      X_1 &&
      X_2 &&
      X_3 &&
      X_4 &&
      X_{N-1} &&
      X_N \\
    };
    \draw (m-2-1) -- (m-1-2);
    \draw (m-2-3) -- (m-1-2);
    \draw (m-2-3) -- (m-1-4);
    \draw (m-2-5) -- (m-1-4);
    \draw (m-2-5) -- (m-1-6);
    \draw (m-2-7) -- (m-1-6);
    \draw (m-2-7) -- (m-1-8);
    \draw (m-2-9) -- (m-1-8);
    \draw (m-2-9) -- (m-1-10);
    \draw (m-2-11) -- (m-1-10);
  \end{tikzpicture}
  \]

  Features that are local to any one of the point clouds will not persist along the zigzag, while global features will be carried along the zigzag to long barcodes. 

  This approach was further studied for practical aspects by \textcite{tausz_applications_2011}, who give concrete algorithms for the computation of the union zigzag, and demonstrate the computational behaviour on a number of concrete examples, including the images dataset studied in \cite{imagespace2008}.

  The union zigzag was also further applied to dynamic network analysis by \textcite{gamble_applied_2012}.
\item[Levelset zigzag] Given a space $\X$ and a continuous function $f:\X\to\R$, the levelset zigzag would relate the levelsets of $f$ through a zigzag, introduced by \textcite{carlsson_zigzag_2009}:
  \[
  \begin{tikzpicture} [thick, ->]
    \matrix (m) [matrix of math nodes,row sep=2ex] 
    {
      & f^{-1}([s_1,s_2]) &
      & f^{-1}([s_2,s_3]) &
      & \dots &
      & f^{-1}([s_{n-1},s_n]) \\
      f^{-1}( s_1) &&
      f^{-1}( s_2) &&
      f^{-1}( s_3) &&
      f^{-1}( s_{n-1}) &&
      f^{-1}( s_n) \\
    };
    \draw (m-2-1) -- (m-1-2);
    \draw (m-2-3) -- (m-1-2);
    \draw (m-2-3) -- (m-1-4);
    \draw (m-2-5) -- (m-1-4);
    \draw (m-2-5) -- (m-1-6);
    \draw (m-2-7) -- (m-1-6);
    \draw (m-2-7) -- (m-1-8);
    \draw (m-2-9) -- (m-1-8);
  \end{tikzpicture}
  \]
  where $a_j$ are the critical values of $f$, and $s_j$ are picked to satisfy:
  \[
  -\infty < s_0 < a_1 < s_1 < a_2 < \dots < s_{n-1} < a_n < s_n < \infty
  \]

  This zigzag produces a computational approach to the interval persistence introduced by \textcite{dey_stability_2007}.
\end{description}

\textcite{carlsson_zigzag_2009} also elaborate the diamond principle to connect it with the Mayer-Vietoris long exact sequence, and give a concrete graphical language for modifying barcodes between union and intersection zigzag sequences. This Mayer-Vietoris relation produces a large diagram from the levelset zigzag (see description below) introduced in the paper that connects to extended persistence, and admits a stability theorem.

\subsubsection{Circular persistence}
\label{sec:circular-persistence}

In a sequence of preprints, \textcites{burghelea_defining_2010,burghelea_persistence_2011} study what they call persistence for circle valued maps. This treats the question of how to adapt the methods of persistent homology in order to deal with studying maps $f:\X\to S^1$ instead of $f:\X\to\R$. Such maps appear naturally when studying cohomology, a fact also underlying the work by \textcite{morozov_persistent_2011} that we described in Section~\ref{sec:cohomology-duality}.

The authors show that by discretizing the map $f:\X\to S^1$ on its critical points, as is done in the real-valued case too, the resulting diagram of homology groups takes the shape of a cyclic quiver: write $G_{2m}$ for a directed graph with $2m$ vertices whose underlying undirected graph is the cycle $C_{2m}$. Then $G_{2m}$ forms a quiver, whose path algebra has representations of the right shape to describe circular persistence. 

Drawing on results by \textcite{donovan_representation_1974} and by \textcite{nazarova_representations_2007}, demonstrating that these quivers have indecomposables classified by barcode spirals coupled with Jordan cells, \textcite{burghelea_persistence_2011} produces algorithms and methods to both compute these indecomposable descriptions, and to solve numerous Betti number computation problems with the spiral and Jordan cell description of a circular persistent homology module. 

\subsection{Diagrams of vector spaces over order categories}
\label{sec:diagr-vect-spac-1}

The approach that shaped Section~\ref{sec:functions-manifold}, studying persistent homology by studying categories of modules over $(\R,\leq)$ is an approach that leads to fruitful approaches to the filtration-based view as well.

\textcite{bubenik_categorification_2012} have results on categorification of persistent homology, see Section~\ref{sec:categorification}. Their approaches, while focused on finite type diagrams over $(\R,\leq)$ seem to be applicable to more general categories of diagrams of vector spaces. With an adapted notion of interleaving distance, the applicability to $(\N,\leq)$-diagrams and thereby to arbitrary filtrations should be immediate.

A recent preprint by \textcite{vejdemo-johansson_interleaved_2012} works out a slightly weakened form of categorical equivalence that relates tame and lower bounded diagrams over $(\R,\leq)$ to tame diagrams over $(\N,\leq)$, thus providing an approach to comparing on a categorical level the different approaches produced by considering filtrations or by considering functions on a manifold. 

Another approach that fundamentally relies on order categories can be found from \textcite{chacholski_combinatorial_2012}. The authors approach multi-dimensional persistence (see Section~\ref{sec:multi-dimens-pers}) by modeling the persistence modules as diagrams over $(\N^r, \leq)$, where $\leq$ is the partial order on $\N^r$ induced by coefficient-wise comparison. They give a concrete and more importantly \emph{local} algorithm for computing the family of invariants $\xi_j$ described by \textcite{carlsson_theory_2009}, thus approaching a more practical and algorithmic approach to persistent homology. 

One more result that emerges from a diagrams of vector spaces approach, and that is important to mention in this paper is from \textcite{ellis_persistent_2011}. The authors consider five important filtrations of finite groups by normal subgroups -- the lower central series, the lower $p$-central series, the derived series, the upper central series, and the upper $p$-central series. For each of these cases, the group cohomology modules of each group in the filtration combine into a \emph{persistent group cohomology diagram} which works as a group invariant with noteworthy discriminatory strength between groups.

\section{Shapes of theories, future directions}
\label{sec:shap-theor-future}

A dichotomy such as the one we have seen above cries out for a unifying theory -- everyone start out with the same underlying problem, and believe they do approximately the same thing, there should be a way to treat all the algebraic foundations in use as aspects of the same underlying theory. While such a unification is not published as this paper is finalized, there are several ongoing efforts in the community that may well lead towards a unifying theory of persistent homology. 

The following descriptions are speculative in nature, describing ongoing work and possible trends, and is fundamentally based on personal communications with 
Gunnar Carlsson, Justin Curry, Robert Ghrist, David Lipsky, Amit Patel, and Primoz Skraba. 
 
With the plethora of differently shaped theories that we have described in Sections \ref{sec:multi-dimens-pers}, \ref{sec:zig-zag-persistence}, \ref{sec:circular-persistence}, and \ref{sec:diagr-vect-spac-1}, a good unification that helps the field forwards will have to deal with the fact that persistent homology is not done with a uni-directional linear progression of some parameter. Instead, a unification will have to systemize handling of differing \emph{shapes} of the theory. 

Once we can accomodate quiver-based shapes, as in zigzag persistence, alongside both continuous and discrete shapes, as with the difference between $(\R,\leq)$-modules and $\kk Q$-modules for an $A_e$-quiver $Q$, the step is not far to start considering tree-like, graph-like, or arbitrary topological spaces describing the underlying shape of a persistence theory.

The group at the University of Pennsylvania, led by Robert Ghrist, has already been building up interest in using sheaves for engineering applications of topology for a while (see \cite{ghrist_applications_2011} -- preprints by Justin Curry and by Sanjeevi Krishnan on their cosheaf and sheaf work are still pending). While the details are yet to be settled, a sheaf-based approach looks promising both for unifying persistent homology and for providing new techniques for applying algebraic topology.

Amit Patel and Robert MacPherson, with input from Paul Bendich, Frédéric Chazal, Herbert Edelsbrunner, Dmitriy Morozov, and Primoz Skraba, are working on using sheaves of \emph{well groups} as a description of persistent homology. Well groups, introduced Edelsbrunner, Morozov and Bendich \cite{well-edelsbrunner-morozov,well-bendich-edelsbrunner-morozov,well-edelsbrunner-morozov-focm}, quantify acceptable noise in the sublevel set approach to $f:\X\to\R$, but have yet to find a large range of computable situations. 

A sheaf-based approach to persistent homology using a particular topological space to describe the shape of the underlying theory is also part of the immediate research agenda for Mikael Vejdemo-Johansson -- classical persistent homology as tame $(\R,\leq)$-diagrams of vector spaces would be the internal algebraic topology of a particular topos of sheaves, where the underlying topological space describes this particular shape of the theory.

\section{Conclusion}
\label{sec:conclusion}

The field of persistent homology draws from a wide range of particular choices of algebraic foundations to describe very similar processes under a common heading. The choices concretely enable a wide range of valuable results, from improved algorithms and new directions of generalization to stability and a road-map towards enabling statistical inference using persistence.

The choices divide, roughly, into two classes with noticable differences -- and from both directions there are things provable in one formalism that are all but inconceivable in the other formalism: stability results seem to be a very odd family of theorems to prove with a strict adherence to a filtration-based point of view, while the results by \textcite{ellis_persistent_2011} are inconceivable if persistent homology can only be thought of as working with functions on a manifold.

Thus, both classes are important view points that enrich the field. Hopefully, the future will bring a satisfactory unification of the foundational choices, demonstrating that there is a single underlying principle to the field.

\section{Acknowledgments}
\label{sec:acknowledgments}

This research was partially funded by the European Union through the project \textsc{Toposys}, grant \# FP7-ICT-318493-STREP.

\clearpage
\printbibliography
\end{document}